\documentclass[11pt,a4paper]{article}

\usepackage[T1]{fontenc}
\usepackage[utf8]{inputenc}
\usepackage{lmodern}
\usepackage[english]{babel}
\usepackage{amsmath,amssymb,amsthm}
\usepackage{booktabs}
\usepackage{geometry}
\usepackage[hidelinks,pdftitle={Fibonacci Enumeration of Residue Classes via Forbidden Binary Factors},pdfauthor={Mike Winkler},pdfsubject={Enumerative combinatorics and parity-word residue classes},pdfkeywords={Fibonacci numbers, binary words, residue classes, finite automata, parity vectors}]{hyperref}
\usepackage{url}

\newcommand{\N}{{\mathbb N}}
\newcommand{\Z}{{\mathbb Z}}

\newtheorem{theorem}{Theorem}[section]
\newtheorem{proposition}[theorem]{Proposition}
\newtheorem{corollary}[theorem]{Corollary}
\newtheorem{lemma}[theorem]{Lemma}
\theoremstyle{definition}
\newtheorem{definition}[theorem]{Definition}
\newtheorem{remark}[theorem]{Remark}

\numberwithin{equation}{section}
\numberwithin{table}{section}

\title{Fibonacci Enumeration of Residue Classes via Forbidden Binary Factors}

\author{Mike Winkler\\[2mm]
	Faculty of Mathematics, Ruhr University Bochum, Germany\\
	mike.winkler@ruhr-uni-bochum.de}

\date{October 1, 2026}

\begin{document}
	\maketitle
	
	\begin{abstract}
		Binary words avoiding adjacent ones are counted by Fibonacci numbers. We identify a natural residue-class enumeration in which this classical word count arises from parity data of the shortcut $3n+1$ map. For an odd starting value, we introduce a local cut time: the first index at which the trajectory enters the residue class $3$ modulo $4$ immediately after an even term. For each finite cut time, the starting values realizing it form a union of complete residue classes modulo a power of two, and the number of these classes is given by an explicit Fibonacci formula. The key step is a recoding of the cut condition as avoidance of the binary factor $11$ in the free part of a parity word. An explicit comparison with the endpoint conventions of the author's 2014 block construction, including the terminal blocks ending at $1$, shows that the resulting residue-class counts prove the two Fibonacci formulas conjectured there in their complete-class form. The construction is a two-state regular-language enumeration; a related theorem of Reyes Jim\'enez has the bitwise-complementary forbidden-factor formulation $00$.
	\end{abstract}
	
	
	
	\section{Introduction}\label{sec:intro}
	
	Binary words avoiding the factor $11$ are counted by Fibonacci numbers. The same recurrence appears whenever a combinatorial object can be encoded by a finite word whose free positions are subject only to this local prohibition. The purpose of this paper is to identify a natural arithmetic instance of this mechanism. A classical parity correspondence converts finite binary words into residue classes modulo powers of $2$, while a new local cut condition singles out precisely those parity words whose free part avoids $11$. The resulting residue-class counts are therefore Fibonacci numbers, and the same viewpoint extends to arbitrary finite-state restrictions through transfer matrices.
	
	The arithmetic source of the construction is the standard shortcut $3n+1$ map
	\[
	T(n)=
	\begin{cases}
		n/2 & \text{if $n$ is even},\\[2pt]
		(3n+1)/2 & \text{if $n$ is odd}.
	\end{cases}
	\]
	For $s\in\N$, write $T^0(s)=s$, $T^k(s)=T(T^{k-1}(s))$, and let $C(s)=(T^k(s))_{k\geq0}$ be the trajectory of $s$. The global $3n+1$ problem asks whether every trajectory reaches $1$; standard references are Lagarias \cite{Lag85,Lag10}. The results below concern only finite parity information and residue-class enumeration. No claim about the global conjecture is used or obtained.
	
	Throughout the paper, $\N=\{1,2,\ldots\}$, and $[a]_m$ denotes the residue class $\{n\in\N:n\equiv a\pmod m\}$. For an odd starting value $s$, let $\ell(s)$ be the first positive index at which the trajectory enters $[3]_4$ immediately after an even term; Definition~\ref{def:ell} gives the formal convention. This single local cut time is the structural device that makes the enumeration possible.
	
	The motivation comes from two Fibonacci formulas conjectured in version 1 (2014) of arXiv:1412.0519 by the author \cite{Win14}. There the formulas were verified computationally for $2\leq h,t\leq50$, but the formulation used two block parameters, two endpoint conventions, and three distinguished residue classes modulo $12$. Version 2 (2021) is a later update of the 2014 manuscript and is not used here. A proof-oriented rewrite appeared as version 3 in August 2026 \cite{Win26v3}; it proved the two return-class counts using parity words and retained the parameters $h$ and $t$, but it did not give a complete two-way comparison with the terminal endpoint conventions of version 1. Because versions 1 and 3 are mathematically different stages of the same arXiv record, they are cited separately throughout. The present revision replaces the two return parameters by the single local cut time $\ell$, states the enumeration first modulo powers of $2$, and then refines it modulo $12$ by the Chinese remainder theorem. Proposition~\ref{prop:bridge2014} supplies the missing comparison, including the blocks ending at $1$.
	
	The parity correspondence used below goes back to Terras \cite{Ter76} and Everett \cite{Eve77}. Fibonacci numbers also occur in other parts of $3n+1$ combinatorics. Albert, Gudmundsson, and \'Ulfarsson \cite{AGU22} studied permutations induced by finite $3n+1$ trajectories; their Fibonacci enumeration holds through length $14$, after which additional permutations occur. More recently, Reyes Jim\'enez \cite{Rey26} proved a Fibonacci theorem for odd starting values whose trajectories avoid the residue class $4$ modulo $6$ over a prescribed interval. His proof uses a directed modular graph and spectral information. In version 3 of that paper, Remark~4.9 gives the equivalent parity-word condition of avoiding the factor $00$ and explicitly relates it to the $11$-avoidance enumeration in \cite{Win26v3}. Section~\ref{sec:finite-state} records this common finite-state context without treating it as an additional principal enumeration result.
	
	Let $(F_m)_{m\geq0}$ be the Fibonacci sequence defined by $F_0=0$, $F_1=1$, and $F_{m+1}=F_m+F_{m-1}$ for $m\geq1$. The two main enumeration results are as follows.
	
	\begin{theorem}\label{thm:intro1}
		For every $r\geq2$, the set
		\[
		\{s\equiv1\pmod4:\ell(s)=r\}
		\]
		is the union of exactly $F_{r-1}$ residue classes modulo $2^{r+2}$.
	\end{theorem}
	
	\begin{theorem}\label{thm:intro3}
		For every $r\geq3$, the set
		\[
		\{s\equiv3\pmod4:\ell(s)=r\}
		\]
		is the union of exactly $F_r-1$ residue classes modulo $2^{r+2}$.
	\end{theorem}
	
	The proof has two ingredients. The finite parity correspondence identifies binary words of length $k$ with residue classes modulo $2^k$. The recoding lemma in Section~\ref{sec:parity} identifies a term in $[3]_4$ with an occurrence of $11$ in the parity word. The free parts of the resulting block words are therefore counted by binary words without adjacent ones. A fixed residue modulo $3$ is independent of the parity word, so the Chinese remainder theorem transfers the same count to each odd residue class modulo $12$.
	
	\section{Parity correspondence and forbidden factors}\label{sec:parity}
	
	For $s\in\N$ and $k\geq0$, put
	\[
	p_k(s)=T^k(s)\bmod2\in\{0,1\},
	\qquad
	P_k(s)=(p_0(s),\dots,p_{k-1}(s))\in\{0,1\}^k.
	\]
	The vector $P_k(s)$ is the parity vector of $s$ of length $k$.
	
	The following finite parity correspondence appears in the early work of
	Terras \cite{Ter76} and Everett \cite{Eve77}. We include a proof in the
	notation used here.
	
	\begin{lemma}[Terras--Everett]\label{lem:terras}
		For every $k\geq1$, the parity vector $P_k(s)$ depends only on the residue
		of $s$ modulo $2^k$, and the induced map
		\[
		\Z/2^k\Z\longrightarrow\{0,1\}^k,
		\qquad [s]\longmapsto P_k(s),
		\]
		is a bijection.
	\end{lemma}
	
	\begin{proof}
		We use induction on $k$. The assertion is clear for $k=1$. Suppose it holds
		for $k$. The map
		\[
		\Phi:\Z/2^{k+1}\Z\longrightarrow\{0,1\}\times\Z/2^k\Z,
		\qquad [s]\longmapsto(p_0(s),[T(s)]),
		\]
		is well-defined. Indeed, replacing $s$ by $s+2^{k+1}m$ changes $T(s)$ by
		$2^km$ in the even case and by $3\cdot2^km$ in the odd case. The map is
		injective. Suppose that $s_1$ and $s_2$ have the same parity and that
		$T(s_1)\equiv T(s_2)\pmod{2^k}$. If both are even, then
		$s_1-s_2=2(T(s_1)-T(s_2))$ is divisible by $2^{k+1}$. If both are odd, then
		$3(s_1-s_2)=2(T(s_1)-T(s_2))$ is divisible by $2^{k+1}$; since $3$ is
		invertible modulo $2^{k+1}$, the same is true of $s_1-s_2$. Thus the inputs
		are congruent modulo $2^{k+1}$. Since the domain and codomain have the same
		cardinality, $\Phi$ is bijective. The induction
		hypothesis applied to $T(s)$ completes the proof.
	\end{proof}
	
	\begin{lemma}[Recoding lemma]\label{lem:recode}
		Let $n$ be odd. Then
		\begin{align*}
			n\equiv1\pmod4&\iff T(n)\text{ is even},\\
			n\equiv3\pmod4&\iff T(n)\text{ is odd}.
		\end{align*}
		Consequently, for every $s\in\N$ and $k\geq0$,
		\begin{align*}
			T^k(s)\equiv3\pmod4
			&\iff p_k(s)=p_{k+1}(s)=1,\\
			T^k(s)\equiv1\pmod4
			&\iff p_k(s)=1,\quad p_{k+1}(s)=0.
		\end{align*}
	\end{lemma}
	
	\begin{proof}
		If $n=4x+1$, then $T(n)=6x+2$ is even. If $n=4x+3$, then
		$T(n)=6x+5$ is odd.
	\end{proof}
	
	Thus a term in $[3]_4$ corresponds to an occurrence of $11$ in the parity
	word, while a term in $[1]_4$ corresponds to an occurrence of $10$.
	
	\begin{corollary}\label{cor:allones}
		For every $k\geq1$,
		\[
		p_0(s)=\cdots=p_{k-1}(s)=1
		\iff s\equiv-1\pmod{2^k}.
		\]
		In particular, every positive integer has an even term in its trajectory.
	\end{corollary}
	
	\begin{proof}
		The first assertion follows by induction. If it holds for $k$, then
		$p_0(s)=\cdots=p_k(s)=1$ is equivalent to $s$ being odd and
		$T(s)\equiv-1\pmod{2^k}$. Hence
		$3s+1\equiv-2\pmod{2^{k+1}}$, which is equivalent to
		$s\equiv-1\pmod{2^{k+1}}$. For $2^k>s+1$, the congruence
		$s\equiv-1\pmod{2^k}$ is impossible.
	\end{proof}
	
	\begin{lemma}\label{lem:words}
		For $m\geq0$, the number of binary words of length $m$ without the factor
		$11$ is $F_{m+2}$.
	\end{lemma}
	
	\begin{proof}
		Let $g_m$ denote this number. We have $g_0=1$ and $g_1=2$. For $m\geq2$,
		a word counted by $g_m$ ends either in $0$, preceded by a word counted by
		$g_{m-1}$, or in $01$, preceded by a word counted by $g_{m-2}$. Thus
		$g_m=g_{m-1}+g_{m-2}$, and the stated formula follows.
	\end{proof}
	
	\section{Cut times and parity blocks}\label{sec:blocks}
	
	For a starting number in $[3]_4$, the trajectory initially contains a
	finite run of consecutive odd terms. We record where this run ends.
	
	\begin{definition}\label{def:a}
		For $s\equiv3\pmod4$, define
		\[
		a(s)=\min\{k\geq1:T^k(s)\equiv1\pmod4\}.
		\]
	\end{definition}
	
	The value $a(s)$ is finite. Indeed, if no such iterate existed, every odd
	iterate would be congruent to $3$ modulo $4$ and would therefore have an odd
	successor by Lemma~\ref{lem:recode}. The parity word would then consist only
	of ones, contrary to Corollary~\ref{cor:allones}. By the same lemma, the
	initial run of ones in the parity word is
	$p_0(s)=\cdots=p_{a(s)}(s)=1$, followed by $p_{a(s)+1}(s)=0$; thus
	$a(s)+1$ is the length of the initial run.
	
	\begin{definition}[Cut time]\label{def:ell}
		Let $s$ be odd. Define
		\[
		\ell(s)=\min\{k\geq1:T^k(s)\equiv3\pmod4
		\text{ and }T^{k-1}(s)\text{ is even}\}.
		\]
		If the set is empty, put $\ell(s)=\infty$. When $\ell(s)<\infty$, the
		finite sequence
		\[
		B(s)=(T^0(s),\dots,T^{\ell(s)-1}(s))
		\]
		is the parity block of $s$. The term $T^{\ell(s)}(s)$ is its cut term.
	\end{definition}
	
	\begin{lemma}[First-cut characterization]\label{lem:firstcut}
		Let $s$ be odd, with the convention that the minimum of an empty set is
		infinity.
		\begin{enumerate}
			\item If $s\equiv1\pmod4$, then
			\[
			\ell(s)=\min\{k\geq1:T^k(s)\equiv3\pmod4\}.
			\]
			\item If $s\equiv3\pmod4$, then
			\[
			\ell(s)=\min\{k>a(s):T^k(s)\equiv3\pmod4\}.
			\]
		\end{enumerate}
	\end{lemma}
	
	\begin{proof}
		If the relevant set on the right is empty, then no cut can occur, so both
		sides are infinite. Suppose first that $s\equiv1\pmod4$, and let
		$k\geq1$ be the first index for which $T^k(s)\equiv3\pmod4$. If $T^{k-1}(s)$ were odd, then
		Lemma~\ref{lem:recode} would give $T^{k-1}(s)\equiv3\pmod4$. For $k=1$
		this contradicts $s\equiv1\pmod4$, and for $k>1$ it contradicts the
		minimality of $k$. Hence $T^{k-1}(s)$ is even, so $k=\ell(s)$.
		
		Now suppose that $s\equiv3\pmod4$, and let $k>a(s)$ be the first later
		index for which $T^k(s)\equiv3\pmod4$. If $T^{k-1}(s)$ were odd, then
		Lemma~\ref{lem:recode} would give $T^{k-1}(s)\equiv3\pmod4$. Since
		$T^{a(s)}(s)\equiv1\pmod4$, we have $k-1>a(s)$, contrary to the
		minimality of $k$. Thus $T^{k-1}(s)$ is even, and again $k=\ell(s)$.
	\end{proof}
	
	This convention is uniform: the cut term is not part of the block. By
	Lemma~\ref{lem:firstcut}, a finite cut time satisfies $\ell(s)\geq2$ when
	$s\equiv1\pmod4$, since $T(s)$ is even. When $s\equiv3\pmod4$, the first
	possible cut occurs after $T^{a(s)}(s)\equiv1\pmod4$ and its even
	successor, so $\ell(s)\geq a(s)+2\geq3$. The ranges in
	Theorems~\ref{thm:intro1} and~\ref{thm:intro3} therefore include every
	possible finite value.
	
	The relation between the cut-time convention and the two endpoint
	conventions of \cite{Win14} is established below in
	Proposition~\ref{prop:bridge2014}.
	
	The cut time $\ell(s)$ need not be finite. The exceptional set is in fact
	nonempty: for
	\[
	s_k=\frac{4^k-1}{3}=1+4+\cdots+4^{k-1}\qquad(k\geq1),
	\]
	we have $T(s_k)=2^{2k-1}$, and hence $\ell(s_k)=\infty$. Nevertheless, every
	exceptional starting value reaches $1$.
	
	\begin{lemma}[Exceptional trajectories]\label{lem:exceptional}
		If $s$ is odd and $\ell(s)=\infty$, then the trajectory of $s$ reaches $1$.
	\end{lemma}
	
	\begin{proof}
		By Lemma~\ref{lem:firstcut}, no later term belongs to $[3]_4$ when
		$s\equiv1\pmod4$, and no term after the initial odd run belongs to
		$[3]_4$ when $s\equiv3\pmod4$. Every subsequent odd term is therefore
		congruent to $1$ modulo $4$. If $o>1$ is such a term, its next odd successor is
		\[
		\frac{3o+1}{2^j}
		\]
		for some $j\geq2$, and is therefore at most $(3o+1)/4<o$. Every positive
		term has a later odd term, since repeated halving of a positive even integer
		reaches an odd integer. A strictly decreasing sequence of positive odd
		integers cannot be infinite. Hence the descent terminates at $o=1$.
	\end{proof}
	
	\begin{definition}[2014 block conventions]\label{def:2014}
		We use the endpoint terminology of \cite[Secs.~4.1--4.5]{Win14} as follows.
		For $s\equiv9\pmod{12}$, the block
		\[
		C^h(s)=(T^0(s),\dots,T^h(s))
		\]
		ends at the least index $h\geq1$ for which
		$T^h(s)\equiv3\pmod4$, provided the trajectory has not reached $1$
		earlier (type $A$); otherwise it ends at the first occurrence
		$T^h(s)=1$ (type $B$). For $s\equiv3,7\pmod{12}$, the block
		\[
		C^t(s)=(T^0(s),\dots,T^t(s))
		\]
		ends at the first even term $T^t(s)\equiv6\pmod8$ that occurs before the
		trajectory reaches $1$ (type $A$), or at the first occurrence
		$T^t(s)=1$ (type $B$). In each case the parameter is the endpoint index.
		
		A residue class modulo $12\cdot2^h$ contained in $[9]_{12}$ is called
		\emph{2014-admissible for $h$} if every starting value in the class has
		2014 parameter $h$. Likewise, a residue class modulo $12\cdot2^{t+1}$
		contained in $[3]_{12}\cup[7]_{12}$ is called \emph{2014-admissible for
			$t$} if every starting value in the class has 2014 parameter $t$. The 2014
		text states that starting values in the listed residue classes have the same
		block length and then counts such admissible classes. That constancy holds
		for the type-$A$ classes; terminal type-$B$ values require separate treatment
		and are not constant on such classes, as Remark~\ref{rem:terminal-representatives}
		shows. We therefore make the complete-class interpretation explicit here; it
		is the interpretation consistent with the tables in \cite[Sec.~4.4]{Win14}.
	\end{definition}
	
	\begin{proposition}[Relation with the 2014 block conventions]
		\label{prop:bridge2014}
		For a starting value whose trajectory reaches $1$, put
		\[
		\eta(s)=\min\{k\geq0:T^k(s)=1\}.
		\]
		Then the following statements hold.
		\begin{enumerate}
			\item If $s\equiv9\pmod{12}$, then the 2014 block is of type $A$ if and
			only if $\ell(s)<\infty$. In that case its parameter is
			$h=\ell(s)$. It is of type $B$ if and only if $\ell(s)=\infty$, and then
			$h=\eta(s)$.
			\item If $s\equiv3,7\pmod{12}$, then the 2014 block is of type $A$ if and
			only if $\ell(s)<\infty$. In that case its parameter is
			$t=\ell(s)-1$. It is of type $B$ if and only if $\ell(s)=\infty$, and then
			$t=\eta(s)$.
			\item The 2014-admissible classes for $h$ are exactly the finite-cut classes
			of Corollary~\ref{cor:mod12} with $c=9$ and $r=h$. The 2014-admissible
			classes for $t$ are exactly the finite-cut classes with
			$c\in\{3,7\}$ and $r=t+1$.
		\end{enumerate}
	\end{proposition}
	
	\begin{proof}
		Once the trajectory reaches $1$, it remains in the cycle $(1,2)$ and no
		later term is congruent to $3$ modulo $4$. Hence a finite cut cannot occur
		after the first occurrence of $1$. Conversely, by
		Lemma~\ref{lem:exceptional}, if $\ell(s)=\infty$ then the trajectory reaches
		$1$. Thus the terminal type $B$ is equivalent to $\ell(s)=\infty$, and a
		finite cut is necessarily of type $A$.
		
		Suppose first that $s\equiv9\pmod{12}$. Then $s\equiv1\pmod4$. By
		Lemma~\ref{lem:firstcut}, when $\ell(s)=r<\infty$ the index $r$ is the
		first $k\geq1$ for which $T^k(s)\equiv3\pmod4$. This is exactly the type
		$A$ endpoint in Definition~\ref{def:2014}, so $h=r=\ell(s)$. If
		$\ell(s)=\infty$, the block ends at the first occurrence of $1$, and hence
		$h=\eta(s)$.
		
		Now let $s\equiv3,7\pmod{12}$ and suppose that $\ell(s)=r<\infty$. By the
		definition of a cut, $T^{r-1}(s)$ is even and $T^r(s)\equiv3\pmod4$.
		Therefore
		\[
		T^{r-1}(s)=2T^r(s)\equiv6\pmod8.
		\]
		No earlier even term can be congruent to $6$ modulo $8$, because such a
		term would have a successor congruent to $3$ modulo $4$ and would therefore
		produce an earlier cut. Hence $T^{r-1}(s)$ is precisely the type $A$
		endpoint of the 2014 block, and $t=r-1=\ell(s)-1$. If
		$\ell(s)=\infty$, the block is terminal and $t=\eta(s)$.
		
		It remains to identify the complete residue classes. First consider the
		$h$-case. Every finite-cut class with $c=9$ and $r=h$ in
		Corollary~\ref{cor:mod12} consists entirely of starting values with
		$\ell(s)=h$; by the first part of the proposition, every element therefore
		has 2014 parameter $h$. Hence each such class is 2014-admissible.
		
		Conversely, suppose a class modulo $12\cdot2^h$ were 2014-admissible for
		$h$ and contained a type-$B$ starting value $s$. Then $T^h(s)=1$. Let $w$
		be the number of odd symbols among the first $h$ parity symbols. Along a
		fixed parity word, composition of the branches of $T$ has the affine form
		\[
		T^h(n)=\frac{3^w n+b}{2^h}
		\]
		with $b$ depending only on that word. Since $s$ and
		$s'=s+12\cdot2^h$ have the same parity word of length $h$,
		\[
		T^h(s')=T^h(s)+12\cdot3^w=1+12\cdot3^w.
		\]
		This value is greater than $1$ and congruent to $1$ modulo $4$. Thus $s'$
		has neither a type-$A$ endpoint nor a type-$B$ endpoint at index $h$, so
		its 2014 parameter is not $h$, contradicting admissibility. Hence an
		admissible $h$-class contains no type-$B$ value. Every element $n$ is therefore
		of type $A$, so $\ell(n)=h$ throughout the class; the class is one of the
		finite-cut classes of Corollary~\ref{cor:mod12}.
		
		The $t$-case is identical. If an admissible class modulo
		$12\cdot2^{t+1}$ contained a type-$B$ value $s$, then $T^t(s)=1$. For
		$s'=s+12\cdot2^{t+1}=s+24\cdot2^t$, the first $t$ parity symbols are
		unchanged and
		\[
		T^t(s')=1+24\cdot3^w\equiv1\pmod8.
		\]
		This is greater than $1$ and is neither a type-$A$ endpoint
		($6\pmod8$) nor the type-$B$ endpoint $1$. Hence the class is not
		admissible, a contradiction. Thus every admissible $t$-class consists of
		type-$A$ values $n$, for which $\ell(n)=t+1$, and is exactly one of the
		finite-cut classes with $c\in\{3,7\}$ and $r=t+1$.
	\end{proof}
	
	\begin{remark}[Terminal representatives]\label{rem:terminal-representatives}
		The complete-class interpretation is essential. Individual type-$B$ starting
		values do occur. If such a value has terminal parameter $k$, then along its
		first $k$ steps
		\[
		1=T^k(s)=\frac{3^w s+b}{2^k}
		\]
		with $w\geq1$ and $b>0$, so $s<2^k/3$. Thus only finitely many terminal
		starting values occur for each fixed parameter.
		
		If $\mathcal U_h$ denotes the union of the 2014-admissible $h$-classes and
		$E_h^B$ the finite set of type-$B$ starting values with terminal parameter
		$h$, then
		\[
		\{s\in[9]_{12}:\text{$s$ has 2014 parameter $h$}\}
		=\mathcal U_h\,\dot\cup\,E_h^B,
		\qquad E_h^B\subset[1,2^h/3).
		\]
		Likewise, with the analogous notation in $[3]_{12}\cup[7]_{12}$,
		\[
		\{s:\text{$s$ has 2014 parameter $t$}\}
		=\mathcal U_t\,\dot\cup\,E_t^B,
		\qquad E_t^B\subset[1,2^t/3).
		\]
		For example,
		\[
		C_B^h(21)=(21,32,16,8,4,2,1)
		\]
		has $h=6$, but $21$ is not part of any 2014-admissible class modulo
		$12\cdot2^6=768$. Indeed, the preceding proof gives
		$T^6(21+768)=37\equiv1\pmod4$, so the parameter is not constant on the
		class $[21]_{768}$. Similarly,
		$C_B^t(3)=(3,5,8,4,2,1)$ has $t=5$, but the class $[3]_{768}$ is not
		2014-admissible for $t=5$. Thus a count of individual representatives would
		not equal the Fibonacci formulas. The conjectures of \cite{Win14} are
		resolved here in the complete-class sense made explicit in
		Definition~\ref{def:2014}.
	\end{remark}
	
	\begin{remark}\label{rem:local}
		Parity blocks are local orbit segments. Distinct blocks may overlap after
		two trajectories merge, and a trajectory may enter a block after its first
		term. These facts do not affect the residue-class enumeration below, which
		depends only on the parity word from the chosen starting number to its
		cut term. A global segment decomposition is therefore unnecessary for
		the present paper.
	\end{remark}
	
	\section{Fibonacci enumeration}\label{sec:counting}
	
	By Lemma~\ref{lem:recode},
	\begin{equation}\label{eq:start}
		\begin{gathered}
			s\equiv1\pmod4
			\iff p_0(s)=1,\quad p_1(s)=0,\\
			s\equiv3\pmod4
			\iff p_0(s)=p_1(s)=1.
		\end{gathered}
	\end{equation}
	
	\begin{proof}[Proof of Theorem~\ref{thm:intro1}]
		By Lemmas~\ref{lem:firstcut} and~\ref{lem:recode}, the condition
		$\ell(s)=r$ is equivalent to
		\begin{equation}\label{eq:ell1}
			(p_k(s),p_{k+1}(s))\neq(1,1)\quad(1\leq k\leq r-1),
			\qquad p_r(s)=p_{r+1}(s)=1.
		\end{equation}
		Together with \eqref{eq:start}, the parity conditions determine a word
		$(p_0,\dots,p_{r+1})$ of length $r+2$. By
		Lemma~\ref{lem:terras}, each admissible word determines one residue class
		modulo $2^{r+2}$.
		
		We have $p_0=1$, $p_1=0$, and $p_r=p_{r+1}=1$. For $r=2$ there is exactly
		one admissible word. If $r\geq3$, condition \eqref{eq:ell1} forces
		$p_{r-1}=0$, and the free letters $p_2,\dots,p_{r-2}$ form a word of
		length $r-3$ (possibly zero) without the factor $11$. Lemma~\ref{lem:words}
		therefore gives $F_{r-1}$ choices, including $F_2=1$ when $r=3$.
	\end{proof}
	
	\begin{proof}[Proof of Theorem~\ref{thm:intro3}]
		Put $a=a(s)$. By Lemmas~\ref{lem:firstcut} and~\ref{lem:recode}, the
		conditions $a(s)=a$ and $\ell(s)=r$ are equivalent to
		\begin{equation}\label{eq:ell3}
			\begin{gathered}
				p_0(s)=\cdots=p_a(s)=1,\qquad p_{a+1}(s)=0,\\
				(p_k(s),p_{k+1}(s))\neq(1,1)
				\quad(a+1\leq k\leq r-1),\\
				p_r(s)=p_{r+1}(s)=1.
			\end{gathered}
		\end{equation}
		Here $a\geq1$ because $s\equiv3\pmod4$, and
		$p_{a+1}=0$ together with $p_r=1$ forces $r\geq a+2$; hence
		$1\leq a\leq r-2$. The value of $a$ is the length of the initial run of
		ones minus one, so it is determined by the word. By
		Lemma~\ref{lem:terras}, every admissible word of length $r+2$ gives one
		residue class modulo $2^{r+2}$.
		
		Fix $a$. Since $p_r=1$, condition \eqref{eq:ell3} for $k=r-1$ forces
		$p_{r-1}=0$. If $r=a+2$ or $r=a+3$, no free letter remains and there is
		exactly one word; note that $F_{r-a-1}=1$ in these cases. If
		$r\geq a+4$, the free letters are $p_{a+2},\dots,p_{r-2}$, a word of
		length $r-a-3$ without the factor $11$. Their number is $F_{r-a-1}$ by
		Lemma~\ref{lem:words}. Since $a$ is determined by the initial run, the
		families for different values of $a$ are disjoint and the counts add. Thus
		\[
		\sum_{a=1}^{r-2}F_{r-a-1}
		=
		\sum_{j=1}^{r-2}F_j
		=F_r-1.
		\]
	\end{proof}
	
	\noindent\textit{Disjointness across cut times.}
	Since $\ell$ is a function, the families of residue classes obtained in both
	Theorems~\ref{thm:intro1} and~\ref{thm:intro3} for distinct values of $r$
	are disjoint. This fact is used in Corollary~\ref{cor:density}.
	
	\begin{corollary}[Refinement modulo $12$]\label{cor:mod12}
		Let $c$ be an odd residue class modulo $12$.
		\begin{enumerate}
			\item If $c\equiv1\pmod4$, then for every $r\geq2$ the set
			\[
			\{s\equiv c\pmod{12}:\ell(s)=r\}
			\]
			is the union of exactly $F_{r-1}$ residue classes modulo
			$12\cdot2^r$.
			\item If $c\equiv3\pmod4$, then for every $r\geq3$ the same set is the
			union of exactly $F_r-1$ residue classes modulo $12\cdot2^r$.
		\end{enumerate}
	\end{corollary}
	
	\begin{proof}
		A class modulo $2^{r+2}$ occurring in Theorem~\ref{thm:intro1} or
		Theorem~\ref{thm:intro3} already fixes the residue modulo $4$. Fixing
		$s\equiv c\pmod{12}$ additionally fixes one residue modulo $3$. Since
		$3$ and $2^{r+2}$ are coprime, the Chinese remainder theorem refines each
		such class to exactly one residue class modulo
		$3\cdot2^{r+2}=12\cdot2^r$. The number of classes is unchanged.
	\end{proof}
	
	\begin{corollary}[The 2014 Fibonacci conjectures in complete-class form]\label{cor:old}
		For every $h\geq2$, the number of 2014-admissible residue classes modulo
		$12\cdot2^h$ in $[9]_{12}$ is $F_{h-1}$. For every $t\geq2$, the number
		of 2014-admissible residue classes modulo $12\cdot2^{t+1}$ in
		$[3]_{12}\cup[7]_{12}$ is $2F_{t+1}-2$.
	\end{corollary}
	
	\begin{proof}
		By Proposition~\ref{prop:bridge2014}, the 2014-admissible $h$-classes are
		precisely the classes with $\ell(s)=h$ in $[9]_{12}$, while the
		2014-admissible $t$-classes are precisely the classes with $\ell(s)=t+1$ in
		$[3]_{12}\cup[7]_{12}$. Corollary~\ref{cor:mod12} therefore gives
		$F_{h-1}$ classes in the first case and
		$2(F_{t+1}-1)=2F_{t+1}-2$ in the second.
	\end{proof}
	
	Table~\ref{tab:classes} lists the first classes in the original parameters
	$h$ and $t$. They agree with the classes tabulated in \cite{Win14}. The first
	count is OEIS A000045 and the second is OEIS A019274 \cite{OEIS}.
	
	\begin{table}[ht]
		\caption{The first 2014-admissible residue classes; the modulus is shown
			in parentheses.}
		\label{tab:classes}
		\centering
		\begin{tabular}{clc}
			\hline
			$h$ & classes in $[9]_{12}$ modulo $12\cdot2^h$ & $F_{h-1}$\\
			\hline
			2 & $9\ (48)$ & 1\\
			3 & $93\ (96)$ & 1\\
			4 & $33,\ 165\ (192)$ & 2\\
			5 & $81,\ 117,\ 237\ (384)$ & 3\\
			6 & $129,\ 333,\ 405,\ 561,\ 645\ (768)$ & 5\\
			\hline
			$t$ & classes in $[3]_{12}\cup[7]_{12}$ modulo $12\cdot2^{t+1}$
			& $2F_{t+1}-2$\\
			\hline
			2 & $27,\ 91\ (96)$ & 2\\
			3 & $19,\ 39,\ 103,\ 147\ (192)$ & 4\\
			4 & $55,\ 67,\ 111,\ 183,\ 195,\ 235,\ 363,\ 367\ (384)$ & 8\\
			\hline
		\end{tabular}
	\end{table}
	
	\begin{corollary}\label{cor:density}
		We have
		\[
		\sum_{r\geq2}\frac{F_{r-1}}{2^r}=1,
		\qquad
		\sum_{r\geq3}\frac{F_r-1}{2^r}=1.
		\]
		Consequently, for every odd residue class $c$ modulo $12$, the
		finite-cut classes have full relative density in $[c]_{12}$; equivalently,
		their complement in $[c]_{12}$ has natural density zero. Every exceptional
		starting number reaches $1$.
	\end{corollary}
	
	\begin{proof}
		The generating function
		\[
		\sum_{m\geq1}F_mx^m=\frac{x}{1-x-x^2}
		\]
		gives $\sum_{m\geq1}F_m2^{-m}=2$. Hence
		\[
		\sum_{r\geq2}\frac{F_{r-1}}{2^r}=1
		\]
		and
		\[
		\sum_{r\geq3}\frac{F_r-1}{2^r}
		=
		\sum_{r\geq3}\frac{F_r}{2^r}
		-
		\sum_{r\geq3}\frac{1}{2^r}
		=
		\left(2-\frac12-\frac14\right)-\frac14=1.
		\]
		Fix an odd residue class $c$ modulo $12$. The class $[c]_{12}$ has density
		$1/12$, while each class modulo $12\cdot2^r$ has density
		$(12\cdot2^r)^{-1}$. If $c\equiv1\pmod4$, then for fixed $R$ the complement
		in $[c]_{12}$ of the classes with $2\leq\ell(s)\leq R$ has density
		\[
		\frac1{12}\left(1-\sum_{r=2}^{R}\frac{F_{r-1}}{2^r}\right),
		\]
		which tends to zero. If $c\equiv3\pmod4$, the analogous density is
		\[
		\frac1{12}\left(1-\sum_{r=3}^{R}\frac{F_r-1}{2^r}\right),
		\]
		which also tends to zero. The exceptional set
		$E_c=\{s\in[c]_{12}:\ell(s)=\infty\}$ is contained in each of these finite
		remainders, according to the value of $c\pmod4$. Hence $E_c$ has upper
		natural density zero, and therefore natural density zero. The final assertion
		is Lemma~\ref{lem:exceptional}.
	\end{proof}
	
	\begin{remark}\label{rem:mean}
		The relative-density weights in Corollary~\ref{cor:density} define
		probability distributions on the finite cut times. Differentiating the same
		generating function gives
		\[
		\sum_{r\geq2}\frac{rF_{r-1}}{2^r}=6,
		\qquad
		\sum_{r\geq3}\frac{r(F_r-1)}{2^r}=8.
		\]
		Thus these distributions have means $6$ on residue classes congruent to
		$1$ modulo $4$ and $8$ on residue classes congruent to $3$ modulo $4$.
		Their masses decay exponentially at rate $\varphi/2$, where
		$\varphi=(1+\sqrt5)/2$.
	\end{remark}
	
	\section{Finite-state context and related work}\label{sec:finite-state}
	
	The preceding proofs are a two-state instance of standard finite-automaton
	counting. If a deterministic finite automaton over $\{0,1\}$ has transfer
	matrix $M$, initial-state vector $u$, and accepting-state vector $v$, then the
	number $a_k$ of accepted words of length $k$ is
	\[
	a_k=u^{\mathsf T}M^k v.
	\]
	By Cayley--Hamilton, $(a_k)$ satisfies the linear recurrence associated with
	the characteristic polynomial of $M$. Combined with
	Lemma~\ref{lem:terras}, the same number counts residue classes modulo $2^k$
	whose parity words satisfy the automaton condition. We use this standard fact
	only to place the preceding enumeration in context.
	
	For words avoiding $11$, two states record whether the current word ends in
	$0$ or in $1$. The transfer matrix is
	\[
	M=\begin{pmatrix}1&1\\1&0\end{pmatrix},
	\]
	with characteristic polynomial $x^2-x-1$. Thus the Fibonacci recurrence in
	Lemma~\ref{lem:words} is exactly the recurrence of this two-state system.
	Theorems~\ref{thm:intro1} and~\ref{thm:intro3} impose fixed prefix and suffix
	conditions; in the second theorem one additionally sums over the initial run
	length.
	
	Reyes Jim\'enez \cite{Rey26} obtains a different Fibonacci enumeration from a
	directed graph of residues modulo $6$. In version 3 of that paper,
	Remark~4.9 shows that, for odd starting values, avoiding the residue class
	$4$ modulo $6$ over a fixed horizon of $m$ parity symbols is equivalent to
	avoiding the factor $00$; among the words of length $m$ beginning in $1$,
	the resulting count is $F_{m+1}$. That remark explicitly relates the $00$
	condition to the $11$-avoidance enumeration in the earlier version of the
	present work \cite{Win26v3}; the two forbidden-factor languages are governed
	by the same Fibonacci transfer matrix. The counted objects are nevertheless
	different: \cite{Rey26} gives a fixed-horizon avoidance count, whereas
	Theorems~\ref{thm:intro1} and~\ref{thm:intro3} are first-passage counts
	with fixed endpoint conditions, and the second theorem also sums over the
	initial run length. The modular-graph method of \cite{Rey26} further yields
	spectral and cycle information not contained in the parity-word count.
	
	\section{Concluding remarks}\label{sec:conclusion}
	
	Theorems~\ref{thm:intro1} and~\ref{thm:intro3} give more than numerical
	identities. They give bijective descriptions of the block classes in terms of
	binary words without adjacent ones. In Theorem~\ref{thm:intro1}, the free
	word has length $r-3$. In Theorem~\ref{thm:intro3}, one first chooses the
	length of the initial run and then a free word without adjacent ones. The
	Fibonacci recurrence is therefore intrinsic to the local block structure.
	The role of the new cut time $\ell$ is essential here: it is the device that
	turns the former block conditions into one uniform finite-state condition.
	Fixing a residue modulo $3$ merely refines each class by the Chinese remainder
	theorem and does not alter the word count.
	
	The finite-state discussion in Section~\ref{sec:finite-state} is contextual: the
	principal enumeration results are Theorems~\ref{thm:intro1} and
	\ref{thm:intro3}, together with the bridge in
	Proposition~\ref{prop:bridge2014}. More complicated local orbit conditions
	recognized by finite automata lead similarly to linear recurrences, but no
	additional enumeration theorem is claimed here.
	
	No claim about the global $3n+1$ problem follows from the finite-state
	enumeration.

\end{document}